\documentclass{amsart}
\usepackage{amsfonts}
\usepackage{amsmath}
\usepackage{amssymb}
\usepackage{amsthm}
\usepackage{soul}
\usepackage[autostyle]{csquotes}
\usepackage{principia}
\MakeOuterQuote{"}
\newcommand{\keyw}[1]{\textit{Keywords---} #1}

\title{On Tarski's Undefinability Theorem}
\author{Stephen Boyce}\thanks{This paper has benefited from use of the principia.sty package.}

\newtheorem{Theorem}{Theorem}[section]
\newtheorem{Corollary}[Theorem]{Corollary}
\newtheorem{Hypothesis}[Theorem]{Hypothesis}
\newtheorem{Proposition}[Theorem]{Proposition}

\theoremstyle{definition}
\newtheorem{Definition}[Theorem]{Definition}
\newtheorem{Example}[Theorem]{Example}
\newtheorem{Fact}[Theorem]{Fact}

\numberwithin{equation}{section}

\begin{document}

\begin{abstract}
  This paper shows that Tarski's revised Undefinability Theorem
  obscures a liar paradox affecting various systems in scope.
  Consider for example Tarski's general theory of classes
  (which contains only variables of finite order) as the object system ($O$),
  and a metatheory $M$ which has a transfinite variable ($Tr^{\omega}$)
  such that (under the intended interpretation): $x \in Tr^{\omega}$ holds iff the
  $O$ formula named by $x$ is true. As there are only a
  denumerable number of true $O$ formulae, it is provable
  within $M$, if we assume the axiom of choice, that there exists a class
  ($X^{4}$), named, under the intended interpretation, by a variable of finite order,
  that contains all, and only, the classes corresponding to the G\"odel numbers of the true $O$ formulae.
  If the class assigned to $Tr^{\omega}$ is well defined
  we may conservatively extend $M$ to $M'$ such that: $M'$ contains
  a new constant of finite type ($Tr^{4}$) and proper axioms
  assigning the class $X^{4}$ to $Tr^{4}$. It is easily shown that
  $M'$ exhibits a liar paradox.  The proof holds, with appropriate
  changes, for ZFC and also implies that
  standard proofs that "arithmetic truth is not arithmetic"
  beg the question of whether the standard interpretation of first-order arithmetic
  is well defined.
\end{abstract}

\maketitle 
\begin{msc}
  03F30, 03B10
\end{msc} \\
\keyw{Tarski's Theorem, First-order arithmetic}
\pagenumbering{arabic}
\section{Introduction}
In his initial formulation, Tarski's Undefinability Theorem asserted that his theory of truth
cannot be applied to provide an adequate definition of truth for languages of infinite order (\cite{tarski1936}: \S 5 Theorem I).
His revised version of this Theorem maintains that, roughly speaking, an adequate definition of
truth for an object system $O$ can be attained in such cases if the order of the metatheory
(in which truth for $O$ is defined) is greater
than that of the object language itself. My aim in the following in twofold: firstly to show that Tarski's
definition of truth, aligned with the revised  Undefinability Theorem, fails to avoid
the liar paradox for cases such as the general theory of classes; secondly,
to show that standard proofs that "arithmetic truth is not arithmetic",
beg the question of whether the standard interpretation of first-order arithmetic
is well defined.

For metatheoretical reasoning, broadly conceived, a primitive classical logical framework is assumed. The various
absurd results derived in what follows are taken as confirmation that the standard
metamathematical definition of truth for a formal theory, due to Tarski, fails for the cases under consideration.
In other words, such pseudo-contradictions are meaningless sequences of symbols indicating
that the attempt to define a meaningful formal theory has failed.

Where convenient, I will often make reference to the meaning associated
with an expression of a formal language under an interpretation as a short-hand way of clarifying
the assumed (formal) definition the expression within the theory. This is
is merely a fa\c{c}on de parler which should not be taken literally; indeed, in view
of the overall argument of the paper, it is not even possible to take such
expressions literally for the cases considered below.

The reader who is familiar with the area may wish to skip to Proposition \ref{proposition_tarski_revised}.

\section{Tarski's Theory of Truth for Formal Languages}\label{section_tarskis_theory_truth}
Let's firstly consider a brief sketch of Tarski's theory of truth for formal languages.
Tarski's definition of truth
is restricted to formal languages as he doubts that truth can be appropriately defined for informal or natural languages
such as English (\cite{tarski1936}: 165). For a formal theory of interest,
defined in a language referred to as the object language,
Tarski presents a metalinguistic method of defining a 
truth predicate in a formal metatheory $M$.
The notions of a formal language ($\mathcal{L}$) and
a formal theory of interest in this language ($T_\mathcal{L}$) are defined in parallel,
such that the combined definition for any particular case includes both:
\begin{enumerate}
\item The syntax of the language $\mathcal{L}$: its primitive symbols
  and the (formation) rules that determine what counts as a well-formed expression / formulae and sentences for this language; and
\item A set of (uninterpreted) statements in $\mathcal{L}$ and also a set of inference rules,
  also defined with respect to the uninterpreted formulae of $\mathcal{L}$, such that:
  the statements stipulate the primitive axioms (or, informally speaking "assumed truths") of the theory $T_\mathcal{L}$ and
  the inference rules specify how a set of $\mathcal{L}$ formulae may be derived or inferred,
  on the basis of either logical or subject-matter specific rules, from another, possibly distinct,
  set of $\mathcal{L}$ formulae (\cite{tarski1936}: 166).
\end{enumerate}
As the definitions of $\mathcal{L}$ and $T_\mathcal{L}$ are linked in this way, I will refer, where convenient, to
either the interpretation of a given formal language, or alternatively of a specific formal theory.
It is always however the interpretation of an associated pair, of $\mathcal{L}$ and $T_\mathcal{L}$,
for which an intended interpretation exists, that is under discussion.

As the need arises on several occasions to use names in one language for
various expressions of a second language, in the context of a sometimes complicated hierarchy of languages,
quasi-quotation (\cite{quine1981}: \S6) is employed to obtain the required precision.
(Alternatives, such as the use of regular quotation marks supplemented by verbose circumlocutions, would very likely result
in miscommunication.) To avoid a misunderstanding, it should be noted that throughout this entire paper
corners are only used for quasi-quotation:
\begin{Definition}[Quasi-quotation]
  Corners - " $\ulcorner$" and "$\urcorner$" - are only used hereon for quasi-quotation (\cite{quine1981}: \S6),
  \emph{not} as a function that returns the G\"odel number of the enclosed character sequence.
\end{Definition}

Tarski's development of key semantic notions - of an interpretation, truth under an interpretation etc. -
for a formal theory or language is similarly premised on the understanding that there exists a given,
intended interpretation of the language of this theory (\cite{tarski1936}: 167). Where convenient hereon the focus will sometimes be
referred to as the object theory ($T_\mathcal{L}$). In contemporary parlance this is usually the set of $T_\mathcal{L}$ sentences that
are formally provable, using the theories axioms and inference rules.
Peano arithmetic (\emph{PA}), or the formal theory of the natural numbers in the language of first-order arithmetic
($\mathcal{L}_{A}$), is a typical example of the combined entity that is the focus of interest.

Given the above syntactical
descriptions of  $\mathcal{L}$ and $T_\mathcal{L}$, Tarski's theory of truth for  $\mathcal{L}$ / $T_\mathcal{L}$ is defined
in a formal metatheory \emph{M}, in a metalanguage ($\mathcal{L}_{M}$) that includes:
the object language ($\mathcal{L}_{O}$), a predicate symbol ($Tr$) for a class which contains specialised
names (in the language of $\mathcal{L}_{M}$) of all, and only,
the $\mathcal{L}_{O}$ sentences that are "true" when the "intended meaning" is applied
to the uninterpreted object language $\mathcal{L}_{O}$. $\mathcal{L}_{M}$ of course contains
much else besides this, including a a predicate symbol ($S$) for a class which contains
all, and only, the names of the $\mathcal{L}_{O}$ sentences; but these are the essential
elements of the definition in $M$ of a true $\mathcal{L}_{O}$ sentence (under a given interpretation):
\begin{quote}
  [Convention T] A formally correct definition of the symbol \textquote*{Tr}, formulated in the metalanguage, will be called \emph{an adequate definition of truth}
  if it has the following consequences:
  \par
  \quad ($\alpha$) all sentences which are obtained from the expression  \textquote*{$x \in Tr$ if and only if $p$} by substituting for the symbol $x$ a structural-descriptive
  name of any sentence of the language in question and for the symbol \textquote*{$p$} the expression which forms the translation of this sentence into the metalanguage;
  \par
  \quad ($\beta$) the sentence \textquote*{for any $x$, if $x \in Tr$ then $x \in S$} (in other words \textquote*{$Tr \subseteq S$}). (\cite{tarski1936}: 187-8).
\end{quote}
As can be seen from this definition, the metalanguage $\mathcal{L}_{M}$ must include, for each sentence
in $\mathcal{L}_{O}$, a sentence which is assigned the same meaning under the interpretation of interest.
This requirement is easily met by defining $\mathcal{L}_{M}$ as an extension of $\mathcal{L}_{O}$.

The assignment of a meanting to object language sentences within the metalanguage
is achieved via set theory using the notion of the satisfaction of a sentential function or sentence
at a denumerable, infinite sequence of objects in a domain of interpretation of the object language.
For example, if our object language is a version of first-order arithmetic and we are using the standard interpretation,
the sentential function "$x_2 = x_3$" is satisfied at
the denumerable sequence $f = \{ 5, 3, 3, 7, 9, \ldots \}$ - since the objects $f_2 = 3, f_3 = 3$ in the sequence assigned to the free variables
($x_2, x_3$) of this sentential function are in the relation assigned to this sentential function (equality defined on natural numbers)
under the interpretation. Using this notion, the notion of the satisfaction of a sentential function $\mathcal{B}(x_{i})$ that is
constructed from other sentential functions - $\mathcal{B}(x_{i}) = \lnot\ \mathcal{C}(x_{i})$,  $\mathcal{B}(x_{i}) = \mathcal{C}(x_{i}) \lor \mathcal{D}(x_{i})$,
or $\mathcal{B}(x_{i}) = (\forall x_J)\mathcal{C}(x_{i}, x_{j})$ etc. - may similarly be defined. (Where no confusion should
arise, I will sometimes suppress corners when, strictly speaking, they ought to be used, as for example
in $\mathcal{B}(x_{i}) = \ulcorner \lnot\ \mathcal{C}(x_{i}) \urcorner$ etc.) From whence one obtains a
definition of the notion that a sentence of $T_\mathcal{L}$ is true under an interpretation if it is
satisfied at every (denumerably) infinite
sequence of objects contained in the domain of interpretation  (\cite{tarski1936}: 195 Def. 23).

Whilst my focus is on Tarski's theory of truth, we should note that the formalisation,
in a metamathematical sense, of both the syntax and semantics of the object theory is
presented. Thus it may be formally proven within the metatheory that a given $T_\mathcal{L}$ sentence
is (or is not) a theorem, that it is (or is not) true under the intended interpretation. Metatheorems
concerning the properties of sets of $T_\mathcal{L}$ sentences that are provable and or true
(under a given interpretation) may also be formally proven within \emph{M}.

To turn now Tarski's Undefinability Theorem let's firstly review the Calculus of Classes, a first-order version of
modern set theory. Whilst the initial version of the Theorem concerns the General Theory of Classes, Tarski's presentation
of the Theory is rather brief and makes use of the definition of the Calculus of Classes.

\section{The Calculus of Classes ($T_{CC}$)}\label{section_first_order_set_theory}
For a brief description of the syntax of the object theory $T_{CC}$ an account of the
primitive symbols, formation rules, logical / proper axioms and inference rules is provided.
Throughout this account some informal comments are made in passing regarding the intended interpretation of the theory,
although, from the metamathematical perspective, precise statements on such matters
should be formulated within a given metatheory.

\paragraph{Primitive Symbols}\label{paragraph_primitive_symbols}
This theories language ($\mathcal{L}_{CC}$) includes a denumerable supply of variables
($x{}_{ _{'}}$, $x_{ _{''}}$, $x_{ _{'''}}$, \ldots) and four signs for constants:
"N", "A", "$\prod$" and "I". Under the intended interpretation,
the constants stand for (respectively): the (classical) logical operations of
negation and alternation (with respect to associated, subsequent formulae),
the universal quantification of the formula(e) in scope (with respect to the individuals ranged over by the
associated, subsequent variables), and the set-theoretic relation of inclusion
(with respect to associated, subsequent variables).

Whilst the primitive symbols for $\mathcal{L}_{CC}$ variables are composed exclusively
from the symbol "$x$" followed by a finite sequence of subscripted strokes, it is sometimes
convenient to use abbreviations for these, such as "$x_2$". Under the intended
interpretation, the domain of interpretation contains "classes of individuals" (\cite{tarski1936}: 169);
hence, informally speaking, the variables range over these.

The meanings of the primitive symbols are possibly easier to grasp if
some (well-formed) formulae are discussed. To consider these
let's firstly review $T_{CC}$'s formation rules.
\paragraph{Formation rules}\label{paragraph_formation_rules}
In contemporary parlance, $T_{CC}$'s atomic formulae may be defined as follows:
\begin{Definition}[Atomic formulae]
  The class of $T_{CC}$ atomic formulae ($AF$) is the smallest set such that,
  whenever $x$ and $y$ are two possibly distinct $\mathcal{L}_{CC}$ variables:
  \begin{equation}
    \ulcorner I x y \urcorner \in AF
  \end{equation}
\end{Definition}
As can be seen from the definition of the atomic formulae,
Tarski uses Polish notation in $\mathcal{L}_{CC}$ so that brackets are not required.
Operators, relations, logical connectives and quantifiers - all of which Tarski refers to as functors
- bind the appropriate number of argument(s) immediately following the
occurrence of the relevant symbol.

The class of $T_{CC}$'s formulae may thus be defined as follows (c.f. \cite{godel1931}: 153):
\begin{Definition}[$\mathrm{Form}(x)$]
  The class of $T_{CC}$ formulae $Form$ is the smallest set such that:
  \begin{enumerate}
  \item $\mathrm{Form}$ includes all atomic formulae:
    \begin{flalign*}
      x \in AF \rightarrow x \in Form
    \end{flalign*}
  \item  Where $x$ and $y$ are any two (possibly distinct) members of $\mathrm{Form}$,
    and $v$ is any $\mathcal{L}_{CC}$ variable:
    \begin{flalign*}
      & \ulcorner N x \urcorner \in Form \\
      & \ulcorner A x y \urcorner \in Form \\
      & \ulcorner \prod v x \urcorner \in Form 
    \end{flalign*}
  \end{enumerate}
\end{Definition}

For the definition of $Tr$ in the metatheory, it is useful to
define the class of $T_{CC}$ sentences. Whilst Tarski's precise,
metatheoretical definition (\cite{tarski1936}: 178, Definition 12) can be used for official purposes,
a brief paraphrase is as follows. Firstly,
to borrow some symbolism from G\"odel (\cite{godel1931}: 167, Def 26), we may
informally define a metatheoretical relation, $v\ \mathrm{Fr}\ x$, that holds of $v$ and
$x$ whenever $v$ occurs as a free variable in $x$:
\begin{Definition}[$v\ \mathrm{Fr}\ x$]
  The relation $v\ \mathrm{Fr}\ x$ is the smallest class
  that includes as members an ordered pair $<v, x>$ if and only if (iff):
  \begin{enumerate}
  \item $v$ is a member of the set of $\mathcal{L}_{CC}$ variables;
  \item $x$ is a  member of the set of $T_{CC}$'s formulae ($\mathrm{Form}$);
  \item there is an occurrence of the variable $v$ in the formula $x$ such
    that, this occurrence is \emph{not} within a subformula $z$
    such that: the formula $\ulcorner \prod v z \urcorner$ is in turn a
    subformula of $x$.
  \end{enumerate}
\end{Definition}
Using this definition of $v\ \mathrm{Fr}\ x$, we may informally
define the class of $T_{CC}$'s sentences as follows:
\begin{Definition}[Sentences]
  The class of $T_{CC}$'s sentences $S$ is the smallest set such that:
  $x$ is a member of $S$ whenever $x$ is a member of the class of
  formulae ($\mathrm{Form}$) that does not have any free variables [$(v)\overline{v\ \mathrm{Fr}\ x}$].
\end{Definition}

To clarify the intended interpretation of $\mathcal{L}_{CC}$,
it is helpful to consider $T_{CC}$'s axioms and inference rules.

\paragraph{Axioms and inference rules}\label{paragraph_axioms_inference_rules}

For the statement of $T_{CC}$'s axioms and inference rules
it is useful to define the class of $T_{CC}$ sentential functions and
the relation of two (possibly distinct) $T_{CC}$ formulae, $x$ and $y$, such that
$x$ is the universal quantification of $y$.
Tarski's definitions (\cite{tarski1936}: 177, Definition 10, \cite{tarski1936}: 176, Definition 8)
can be used for official purposes; a paraphrase however is as follows. 
\begin{Definition}[$SFx$]
  The class of $T_{CC}$'s sentential functions $SFx$ is the smallest set such that:
  $x$ is a member of $SFx$ whenever $x$ is a member of the class of
  formulae ($\mathrm{Form}$) that is \emph{not} a member of the class of $T_{CC}$'s sentences $S$.
\end{Definition}
The notion of universal quantification defined may be paraphrased as follows. 
\begin{Definition}[Universal quantification]
  $x$ is the universal quantification of $y$ iff:
  \begin{enumerate}
  \item $y$ is a  member of the set of $T_{CC}$'s formulae ($\mathrm{Form}$); and
  \item $x$ is a  member of the set of $T_{CC}$'s sentences ($S$); and \\
    either:
  \item $x$ is equal to $y$; or \\
    all of the following hold:    
  \item there exists a sequence of $n$ $\mathcal{L}_{CC}$ variables $v_1 , v_2, \ldots v_n$;
  \item there exists a sequence of $n+1$ $\mathcal{L}_{CC}$ formulae $z_1 , z_2, \ldots z_{n+1}$;
  \item the first formula in this sequence of formulae is $y$;
  \item for $1 < i \leq n+1$, the $i^{\text{th}}$ formula is given by:
    \begin{equation}
      z_i =  \ulcorner \prod v_i z_{i-1} \urcorner
    \end{equation}
  \item $x$ is the formula $z_{n+1}$.\
  \end{enumerate}
\end{Definition}
With the above definitions, we may define $T_{CC}$'s logical axioms (\cite{tarski1936}: 179, Definition 13 $\alpha$) as follows:
\begin{Definition}[Logical axioms]\label{definition_cc_logical_axioms}
  $x$ is a logical axiom of $T_{CC}$ iff:
  \begin{enumerate}
  \item $x$ is a  member of the set of $T_{CC}$'s sentences ($S$); and \\
  \item either:
    \begin{enumerate}
    \item there exist $T_{CC}$ sentences $y$, $z$, $u$ such that $x$ is equal to
      one of the following $T_{CC}$ sentences:
      \begin{flalign}
        & \ulcorner A NAyy y \urcorner\\
        & \ulcorner A Ny Ayz \urcorner \\
        & \ulcorner A NAyz Ayz \urcorner\\
        & \ulcorner A  NANyz A NAuy Auz \urcorner
      \end{flalign}
    \item or else there exist $T_{CC}$ sentential functions
      (members of $SFx$) $y$, $z$, $u$ such that $x$ is equal to
      the universal quantification of one of the following $T_{CC}$ sentential functions:
      \begin{flalign}
        & \ulcorner A NAyy y \urcorner\\
        & \ulcorner A Ny Ayz \urcorner \\
        & \ulcorner A NAyz Ayz \urcorner\\
        & \ulcorner A  NANyz A NAuy Auz \urcorner
      \end{flalign}
    \end{enumerate}
  \end{enumerate}
\end{Definition}
As the axioms are adapted from the propositional calculus of \emph{Principia Mathematica},
it may improve the readability of the above definition if corresponding axioms, as stated in \emph{Principia},
are examined (\cite{pm1910v1}: 100-1):
\begin{flalign}
  & \text{Taut}\ \pmast 1\pmcdot2. \quad \pmthm \pmdott p \pmor p  \pmdot \pmimp \pmdot  p \\
  & \text{Add}\ \pmast 1\pmcdot3. \quad \pmthm \pmdott q  \pmdot \pmimp \pmdot  p \pmor q \\
  & \text{Perm}\ \pmast 1\pmcdot4. \quad \pmthm \pmdott p \pmor q \pmdot \pmimp \pmdot  q \pmor p \\
  & \text{Sum}\ \pmast 1\pmcdot6. \quad \pmthm \pmdottt q \pmimp r \pmdot \pmimp \pmdott p \pmor q \pmdot \pmimp \pmdot  p \pmor r 
\end{flalign}
In stating the proper axioms $T_{CC}$ (\cite{tarski1936}: 179, Definition 13 $\beta$), to improve the readability, I will translate the meaning
of the axioms, under the intended interpretation, into English (c.f. \cite{tarski1936}: 170).
\begin{Definition}[Proper axioms]
  $x$ is a proper axiom of $T_{CC}$ iff $x$ is a member of the set of $T_{CC}$'s sentences ($S$) that,
  under the intended interpretation, translates to one of the following propositions (where
  various standard abbreviations have been introduced, $x{}_{ _{'}}$ translates to $a$,
  $x{}_{ _{''}}$ translates to $b$ etc.,
  and the variables range over the sets of individuals in the domain of interpretation):
  \begin{equation}
    (\forall a) (a \subseteq a)
  \end{equation}
  \begin{equation}
    (\forall a) (\forall b) (\forall c) (a \subseteq b \rightarrow b \subseteq c \rightarrow a \subseteq c)
  \end{equation}
  \begin{equation}
    (\forall a) (\forall b) (\exists c) [a \subseteq c \land b \subseteq c \land
      (\forall d) (a \subseteq d \rightarrow b \subseteq d \rightarrow c \subseteq d)] 
  \end{equation}
  \begin{equation}
    (\forall a) (\forall b) (\exists c) [c \subseteq a \land c \subseteq b \land
      (\forall d) (d \subseteq a \rightarrow d \subseteq b \rightarrow d \subseteq c)] 
  \end{equation}
  \begin{flalign*}
  \quad   \{ (\forall a) (\exists b) & \{(\forall c) (\forall d) [(c \subseteq a \rightarrow c \subseteq b \rightarrow c \subseteq d)\ \land \tag{3.15} \\
      & (a \subseteq c \rightarrow b \subseteq c \rightarrow d \subseteq c)]\ \land \\
    & (\forall e) [e \subseteq b \lor (\exists f) (f \subseteq a \land \overline{f \subseteq b} \land f \subseteq e)] \}
  \end{flalign*}
\end{Definition}
The proper axioms, informally assert that:
\begin{enumerate}
\item The subset relation is reflexive, as every set is identical with itself;
\item Among sets of individuals, the subset relation is transitive;
\item For each pair of sets of individuals, there exists a set that contains the union of both sets;
\item For each pair of sets of individuals, there exists a set that contains the intersection of both sets;
\item Any set of individuals is either the empty set, a set that contains all individuals,
  or a set that contains one or more individuals that some other set(s) do not contain.
\end{enumerate}
Tarski refers to the logical and proper axioms of $T_{CC}$ as primitive sentences (\cite{tarski1936}: 179).
This is similar to \emph{Principia}'s category of Primitive propositions (\cite{pm1910v1}: $\pmast 1$), with the qualification
that the axioms, in the former case, are defined with respect to uninterpreted sequences of signs (albeit
with an intended interpretation in mind).

The inference rules of $T_{CC}$ may be summarised as follows.
\begin{Definition}[Inference rules]
  The rules of inference of $T_{CC}$ (\cite{tarski1936}: 181-2, Definition 15 $\gamma$-$\zeta$) provide for:
  \begin{description}
  \item[$\gamma$] inferences involving a change of variables, where a clash of variables does not occur;
  \item[$\delta$] modus ponens or detachment as an inference rule;
  \item[$\epsilon$] distribution of a universal quantifier into a disjunction;
  \item[$\zeta$] exportation of a universal quantifier from inside the scope of
    a disjunction.
  \end{description}
\end{Definition}

Tarski's initial version of his theory of truth proposed that languages
of an infinite order required a different treatment to other categories of
language. To examine this argument let's firstly consider the assumed semantic category
account of the order of a formal language, noting of course that
Tarski revises his method of assigning an order to a formal language
of infinite order (\cite{tarski1936}: 268-9).

\section{The Order of a Language based on its Semantic Category}\label{section_order}
If a language is assigned an order
on the basis of the semantical category of expressions it contains (\cite{tarski1936}: \S4),
then a first-order language is as follows.\footnote{The definition
agrees with contemporary usage e.g. (\cite{mendelson2015}: 53 footnote "*").}
\begin{Definition}[First-order formal language]
  A formal language is first-order iff all of the following conditions are met:
  \begin{enumerate}
  \item All variables of the language are assigned, under under any interpretation, to individuals in the domain of interpretation. 
  \item None of the predicate and function symbols of the language are variables. Under any interpretation, each is assigned a fixed meaning.
    The interpretation of a given predicate (or function) symbol generally associates with the symbol
    a class of $n$-tuples of individuals from the domain of interpretation
    (or, respectively a function that takes $n$-tuples of individuals as arguments) - where,
    for any given predicate or function symbol, $n$ is a fixed natural number greater than or equal to zero).
  \item Apart from the symbols for logical constants (including quantifiers), predicates and functions, the only symbols for constants
    are, under any given interpretation, assigned to individuals in the domain of interpretation.
  \end{enumerate}
\end{Definition}

The definition aims simply to clarify some of the key features of
commonly encountered classical first-order theories.
It is not intended to cover all kinds of formal languages or all complications that
would need to be addressed for a comprehensive definition.
It does not for example deal with the case of a purely syntactical formal language or
formal languages in the case where some non-standard interpretation is applied. (To avoid
circumlocutions referring to an assigned interpretation, we could alternatively
adopt Tarski's approach and assign symbols to categories on the syntactical basis of
whether they are arguments to other symbols in an expression, or the category of symbol
that they take as argument in forming an expression and so on.)

Clause $\alpha$ of Convention T illustrates one kind of constant excluded by the third condition
in the above definition. A formal metalanguage ($\mathcal{L}_{M}$) defining truth for an object language
$\mathcal{L}_{O}$ may include a functor $\mathfrak{Tr}(x, p)$ that,
under the intended interpretation of $\mathcal{L}_{M}$, takes the value true for a pair of arguments iff all of the following three conditions hold
simultaneously: (i) the first argument to the functor is in $Tr$, the set of names of object language ($\mathcal{L}_{O}$) sentences
that are true under the standard interpretation;
(ii) the first argument is also the $\mathcal{L}_{M}$ name for an $\mathcal{L}_{O}$ sentence that is the second argument;
(iii) if the first two conditions hold, this implies that second argument is an $\mathcal{L}_{O}$ sentence
that is true under the intended interpretation of $\mathcal{L}_{O}$.
By Tarski's account, it is necessary, if a liar paradox is to be avoided, that the order of
the metalanguage $\mathcal{L}_{M}$ \emph{must} be strictly greater than that of $\mathcal{L}_{O}$.

As the illustration highlights, in Tarski's account of truth the assignment of an order to a language (and associated formal theory), is
designed to ensure that the definition of a truth predicate does not enable the construction of a sentence that asserts of itself that it is false.

To continue up the chain of orders, a language is, roughly speaking, second order iff it has, in addition to individual variables, variables
that, under the intended interpretation, range over: either (i) the properties and / or relations of individuals in the domain of interpretation; and / or
(ii) functions of individuals in the domain of interpretation. A second-order language may also of course have first-order and second-order predicate and function symbols
that, under the intended interpretation, are assigned a constant meaning (in the sense explained in the definition of a first-order language).
Similarly, a "third order" language has, in addition to the components of a second order language,
variables that range over either the second-order relations and / or second-order functions, and so on.

If the upward progression of types of variables in a language
continues with no finite limit then the order of the language is infinite (\cite{tarski1936}: 220).
(As has just been illustrated, where convenient hereon properties of a given order shall be
taken to be included in relations of that order.)
The language of the General Theory of Classes discussed in the following section is
an example of a language of infinite order in this sense.

Tarski's subsequent revision of this method of assigning an order to a language is concerned primarily with
languages of infinite order. Under the revised approach, transfinite ordinals are used to
assign an order to languages of infinite order, rather assigning all to a single category (\cite{tarski1936}: 268-9).
Before considering this revision, let's consider why Tarski initially took the view that
languages of infinite order are out of scope of his theory of truth.

\section{Tarski's Initial Theorem on languages of infinite order}\label{section_languages_infinite_paradox}
For the case of languages of infinite order that are in scope for G\"odel's proof, Tarski outlines a proof of the following Theorem (\cite{tarski1936}: \S5 Theorem I), which
states in brief that his theory of truth cannot be applied to provide an adequate definition of truth for such languages:
\begin{quote}
  [Theorem I] \quad ($\alpha$) In whatever way the symbol \textquote*{Tr}, denoting a class of expressions, is defined in the metatheory, it will be possible to derive from it
  the negation of one of the sentences which were described in the condition ($\alpha$) of the Convention T;
  \par
  \quad ($\beta$) assuming that the class of all provable sentences of the metatheory is consistent, it is impossible to construct an adequate definition of truth
  in the sense of Convention T on the basis of the metatheory. (\cite{tarski1936}: 247).
\end{quote}
Tarski's proof of this Theorem focuses on the case of the General Theory of Classes,
the syntax of which is quite similar to G\"odel's system $P$ (\cite{godel1931}).
To review Tarski's proof, the following subsection provides a brief
summary of the Theory's syntax.

\subsection{The General Theory of Classes ($T_{GTC}$)}

To briefly review the syntax of $T_{GTC}$ I will focus on the
primitive symbols, formation rules, logical / proper axioms and inference rules.
For future reference, let's firstly introduce some symbolism relating to the formal metatheory
of $T_{GTC}$.
\begin{Definition}[The formal metatheory ${T_M}_{GTC}$ of $T_{GTC}$]
Let ${T_M}_{GTC}$ be the formal metatheory, in the language ${\mathcal{L}_M}_{GTC}$,
in which the syntax and semantics of $T_{GTC}$ - including
the notion of truth under the intended interpretation are formally defined.
\end{Definition}
Returning to $T_{GTC}$ itself, the language $\mathcal{L}_{GTC}$
and syntax more generally of this theory
are very similar to those of $T_{CC}$.
\paragraph{Primitive Symbols}
The language $\mathcal{L}_{GTC}$) includes
the same four signs for (logical) constants
as the calculus of classes: "N", "A", "$\prod$" and "I",
associated with the same meaning under the intended interpretation.
 $\mathcal{L}_{GTC}$) includes a denumerable supply of variables
for each order:
\begin{flalign*}
  X_{'}^{'}, X_{'}^{''}, X_{'}^{'''}, \ldots \\
  X_{''}^{'}, X_{''}^{''}, X_{''}^{'''}, \ldots \\
  X_{'''}^{'}, X_{'''}^{''}, X_{'''}^{'''}, \ldots \\
  \ldots
\end{flalign*}
Whilst the primitive symbols for $\mathcal{L}_{CC}$ variables are composed exclusively
from the symbol "$X$" followed by a finite sequence of subscripted strokes, and
a finite sequence of superscripted strokes, where convenient abbreviations
such as "$X_2^{5}$" etc. will be used. $X_{k}^{i}$ is thus the $i^{\text{th}}$ variable of the
$k^{\text{th}}$ order. The last ellipsis in the above listing of
variables is intended to indicate that for each natural number $k$
there exists a denumerable sequence of $\mathcal{L}_{CC}$ variables
for that order: $X_{k}^{1}, X_{k}^{2}, X_{k}^{3} \ldots$.

$T_{GTC}$'s formation rules are very similar to those of $T_{CC}$.
\paragraph{Formation rules}
$T_{GTC}$'s atomic formulae may be defined thus:
\begin{Definition}[Atomic formulae $T_{GTC}$]
  The class of $T_{GTC}$ atomic formulae ($AF_{GTC}$) is the smallest set such that,
  whenever $x$ and $y$ are two necessarily distinct $\mathcal{L}_{GTC}$ variables such that
  the order of $x$ is one greater than the order of $y$:
  \begin{equation}
    \ulcorner x y \urcorner \in AF_{GTC}
  \end{equation}
\end{Definition}

Corresponding then to definitions provided above for $T_{CC}$, with appropriate changes,
we may then define the following classes / relations for $T_{GTC}$:
\begin{enumerate}
\item $\mathrm{Form}_{GTC}(x)$: the class of $T_{GTC}$'s formulae;
\item $S_{GTC}$: the class of $T_{GTC}$'s sentences;
\item $SFx_{GTC}$: the class of $T_{GTC}$'s sentential functions;
\item $x$ is a universal (or existential) quantification of $y$ ($T_{GTC}$):
  the relation that holds between two (possibly distinct) $T_{GTC}$
  formulae $x$ and $y$ when they are either the same or else
  $x$ is a sentence that is
  the universal (or existential) quantification of $y$
  with respect to one or more $T_{GTC}$ variables
  that have free occurrences in $y$.
\end{enumerate}

$T_{GTC}$'s axioms and inference rules are also similar to $T_{CC}$'s.
$T_{GTC}$'s logical axioms are the $T_{GTC}$ sentences which,
in the appropriate sense, correspond to the axioms of sentential calculus
(possibly followed by universal quantification).
\begin{Definition}[Logical axioms ($T_{GTC}$)]
  $x$ is a logical axiom of $T_{GTC}$ iff: the Definition (\ref{definition_cc_logical_axioms}) of a $T_{GTC}$ Logical axiom
  is satisfied, with "$T_{CC}$" replaced throughout by "$T_{GTC}$'s".
\end{Definition}

\begin{Definition}[Proper axioms ($T_{GTC}$)]
The proper axioms of $T_{GTC}$ are defined, with corners omitted, as follows  (\cite{tarski1936}: 243):
\begin{enumerate}
\item \textbf{Pseudodefinitions (comprehension axioms)} are all sentences of the following form
  (where $\mathcal{B}$ is any sentential function which does not contain the free variable $X_{k}^{n+1}$):
  \begin{flalign*}
    N \prod X_{k}^{n+1} N \prod X_{l}^{n} A\ NA\ NX_{k}^{n+1} X_{l}^{l}\ N\mathcal{B}\
    NA\ NNX_{k}^{n+1} X_{l}^{l}\ NN\mathcal{B}
  \end{flalign*}
 In contemporary notation, where $\mathcal{B}$ again is any sentential function which does not contain the free variable $a^{n+1}$:
    \begin{flalign*}
    (\exists a^{n+1}) (\forall b^{n}) \{[b^{n} \in a^{n+1} \land \mathcal{B}] \lor [\lnot(b^{n} \in a^{n+1}) \land \lnot \mathcal{B}] \}
  \end{flalign*}
\item \textbf{Laws of extensionality (equality axioms)} are all sentences of the following form (asserting that classes with the same elements are equal):
    \begin{flalign*}
    \prod X_{k}^{p+2} \prod X_{l}^{p+1} \prod X_{m}^{p+1} \{A\ \\
    N\prod X_{n}^{p}N [A\ (NA\ NX_{l}^{p+1} X_{n}^{p}\ NNX_{m}^{p+1} X_{n}^{p})\ \\
      (NA\ NNX_{l}^{p+1} X_{n}^{p}\ NX_{m}^{p+1} X_{n}^{p})] \\
    [A\ NX_{k}^{p+2} X_{l}^{p+1}\ X_{k}^{p+2} X_{m}^{p+1}]\}
  \end{flalign*}
  In more contemporary notation:
  \begin{flalign*}
    (\forall a^{n+2}) (\forall b^{n+1}) (\forall c^{n+1})\
    \{& (\exists d^{n}) [(\{d^{n} \in b^{n+1}\} \land \lnot \{d^{n} \in c^{n+1}\})\ \\
      & \lor (\lnot \{d^{n} \in b^{n+1}\} \land \{d^{n} \in c^{n+1}\})] \\
    & \lor\ [(b^{n+1} \in a^{n+2}) \rightarrow (c^{n+1} \in a^{n+2})] \}
  \end{flalign*}
\item A single \textbf{axiom of infinity}, guaranteeing the existence of
  an infinite number of \textit{individuals}. An axiom for Tarski's system, in contemporary
  notation is as follows:\footnote{Tarksi's statement of this axiom is as follows (\cite{tarski1936}: 243):
  \begin{align*}
    \textstyle \bigcup_{1}^{3}( \bigcup_{1}^{2} \epsilon_{1,1}^2\ .\
    \bigcap_{1}^{2}(\overline{\epsilon_{1,1}^2}\ +\
    \bigcap_{2}^{2}(\epsilon_{1,2}^2\ .\
    \bigcap_{1}^{1}(\epsilon_{1,1}^1\ +\ \overline{\epsilon_{2,1}^1})\ .\
    \bigcup_{1}^{1}(\epsilon_{1,1}^1\ +\ \overline{\epsilon_{1,2}^1}))))
  \end{align*}}
  \begin{flalign*}
    (\exists a^{3}) (\exists a^{2}) (a^{2} \in a^{3}\ \land (\exists a^{1})(a^{1} \in a^{2} \land 
    (\forall b^{2}) (b^{2} \in a^{3}\ \rightarrow \\
      (\exists c^{2}) (c^{2} \in a^{3}\ \land
      (\forall a^{1}) (a^{1} \in b^{2} \rightarrow a^{1} \in c^{2})\ \land \\
      (\exists c^{1})(c^{1} \in c^{2} \land\ c^{1} \notin b^{2})))))
  \end{flalign*}
  The axiom states that there's an order-three class ($a^{3}$) which contains
  a non-empty order-two class ($a^{2}$), and for every order-two class ($b^{2}$) in $a^{3}$,
  there exists a (necessarily distinct) order-two class ($c^{2}$) that contains
  every individual in $b^{2}$ and at least one individual not in $b^{2}$.
\end{enumerate}
\end{Definition}

$T_{GTC}$'s inference rules correspond to those of $T_{CC}$, with appropriate changes:
\begin{Definition}[Inference rules ($T_{GTC}$)]
  The rules of inference of $T_{GTC}$ provide for a change of variables
  (where a clash of variables does not occur), modus ponens,
  distribution of a universal quantifier into a disjunction,
  and exportation of a universal quantifier from inside the scope of a disjunction .
\end{Definition}

The above provides enough information about $T_{GTC}$ to review Tarski's proof of
the first Theorem of \S5 (\cite{tarski1936}: 247).

\subsection{A review of Tarski's Retracted Undefinability Theorem (\S5 Th 1)}\label{subsection_retracted_undefineability}

The proof assumes as a hypothesis, in order to derive a contradiction,
that there exists a formal metatheory ${T_M}_{GTC}$ which
defines truth for $T_{GTC}$ in a manner that satisfies the requirements of Convention T.

Throughout the proof, where convenient,
the natural number $n$ is identified with the class of all classes of
individuals that have $n$ members. Where convenient it is also assumed that
a formal \emph{metatheory of} ${T_M}_{GTC}$ - hence, relative
to $T_{GTC}$, a metametatheory - is available. The proof,
if presented as a deduction in a formal theory, would be
carried out in this metametatheory, though to simplify the presentation
it is described as though the formal deductions etc.
are to be performed in the metatheory ${T_M}_{GTC}$.

Let $\phi_0$, $\phi_1$, \ldots, $\phi_n$, \ldots
be an enumeration of the expressions of the object
language ($\mathcal{L}_{GTC}$) via a G\"odel numbering (\cite{tarski1936}: 249).
Let us also assume that this arithmetization of $\mathcal{L}_{GTC}$
is formally defined within our metatheory ${T_M}_{GTC}$, so that
it may be formally proved in ${T_M}_{GTC}$ that the number $j$
is thus associated with the $\mathcal{L}_{GTC}$ expression $\phi_j$ and vice versa.

For $k$ a metametalinguistic variable that ranges over the natural numbers only,
let the metametalinguistic expression $\ulcorner \iota_k \urcorner$
name a metalinguistic (${T_M}_{GTC}$) expression (in ${\mathcal{L}_M}_{GTC}$)
for the corresponding $T_{GTC}$ expression (in $\mathcal{L}_{GTC}$)
such that: 
this $T_{GTC}$ expression, informally speaking, asserts,
under the intended interpretation, that
the class assigned to the first $\mathcal{L}_{GTC}$ order-three variable
$X_{1}^{3}$ names the (natural) number $k$. (This informal gloss can of course be expressed more precisely
- the sentential function in question is satisfied at a (denumerably) infinite sequence
$f$ of classes of classes of individuals (in the domain of interpretation)
if the first element of the sequence is the class of classes of such individuals that have
exactly $k$ members.)

The proof of Theorem I (\S5) focuses on the metametalinguistic expression \\
"$\bigcup_{1}^{3}(\iota_{n} . \phi_n)\ \overline{\in}\ Tr$".
"$\bigcup_{i}^{k}$" and "$\bigcap_{i}^{k}$" are,  respectively,
the ${T_M}_{GTC}$
names for the $T_{GTC}$ expressions for existential
and universal quantification with respect to
the $i^{\text{th}}$ variable of the $k^{\text{th}}$ order,
i.e. $N\prod X_{i}^{k} N$ and $\prod X_{i}^{k}$.

The dot, ``$.$'', is used in the ${T_M}_{GTC}$,
and also in the metametalanguage, to
form names for $T_{GTC}$ formulae that are conjunctions.
So for example where $x$ and $y$ are metalinguistic names for the
expressions "$X_{'}^{'''} X_{'}^{''}$" and "$X_{''}^{''} X_{'}^{'}$" respectively,
$\ulcorner x . y \urcorner$ names the following expression
\begin{equation*}
  NA NX_{'}^{'''} X_{'}^{''} NX_{''}^{''} X_{'}^{'}
\end{equation*}

Putting the above together,
the metametalinguistic expression
"$\bigcup_{1}^{3}(\iota_{n} . \phi_n)\ \overline{\in}\ Tr$"
is thus a sentential function, that, informally speaking, takes names for natural numbers
as arguments. For a particular natural number $k$, the value of the function
is a ${T_M}_{GTC}$ sentence $\mathcal{B}$ such that:
where $\phi_n$ is a $T_{GTC}$ sentential function with "$X_{'}^{'''}$" as it's only
free variable,  $\mathcal{B}$ asserts under the intended interpretation,
that the sentential function $\phi_n$ fails to hold for the number $n$, i.e. it is not
satisfied at any (denumerably) infinite sequence $f$ of classes of classes of individuals
(in the domain of interpretation)
if the first element of the sequence is the class of classes of such individuals that have
exactly $n$ members.

Since however both $\mathcal{L}_{GTC}$ and the language of the metametatheory
include $\mathcal{L}_{GTC}$ the sentential function "$\bigcup_{1}^{3}(\iota_{n} . \phi_n)\ \overline{\in}\ Tr$"
is itself one of the expressions $\phi_k$, for a given choice of $k$,
with "$X_{'}^{'''}$" as it's only free variable.

Moreover as Tarski observes:
\begin{quote}
  \ldots with this function we can correlate another function which is equivalent to it for any value
  of "$n$", but which is expressed completely in terms of arithmetic. We shall write this new
  function in the schematic form "$\psi(n)$" (\cite{tarski1936}: 250, modified: single quotation marks changed to double)
\end{quote}

Tarski's original version of the demonstration provides only an informal argument to support this claim,
though accepted methods of constructing the required function
are generally covered in more contemporary presentations. To illustrate one approach to
filling in some of the omitted details, let us firstly
clarify the requirements for the sought after,
arithmetically defined (meta)metatheoretical sentential function $\psi(n)$.

For a given well-defined class of $T_{GTC}$ sentences that are true under the intended interpretation,
an arithmetic definition of a (meta)metatheoretical sentential function $\psi(n)$
will be (truth functionally) equivalent to $\bigcup_{1}^{3}(\iota_{n} . \phi_n)\ \overline{\in}\ Tr$
provided that:\\
\begin{enumerate}
\item $\psi(n)$ is true, under the intended interpretation, for all, and only, $n \in \mathbb{N}$ such that:
  \begin{enumerate}
  \item $\phi_n$ is an $\mathcal{L}_{GTC}$ sentential function such that:
    ($\alpha$) $\phi_n$ has "$X_{'}^{'''}$" as its only free variable; and 
    ($\beta$) roughly speaking: under the intended interpretation, $\phi_n$
    is \emph{not} satisfied at any (denumerably) infinite sequence
    $f$ of classes of classes of individuals (in the domain of interpretation)
    when the first element of the sequence is the class of classes of such individuals that have
    exactly $n$ members.
  \item $\phi_n$ is a $\mathcal{L}_{GTC}$ sentence that is false
    under the intended interpretation, or, equivalently, is not in $Tr$;
  \end{enumerate}
\item $\psi(n)$ is false, under the intended interpretation, for all, and only, $n \in \mathbb{N}$ such that:
  \begin{enumerate}
  \item $\phi_n$ is neither an $\mathcal{L}_{GTC}$ sentence,
    nor an $\mathcal{L}_{GTC}$ sentential function, or else
    is an $\mathcal{L}_{GTC}$ sentential function that
    does not have "$X_{'}^{'''}$" as its only free variable; 
  \item $\phi_n$ is an $\mathcal{L}_{GTC}$ sentential function such that:
    ($\alpha$) $\phi_n$ has "$X_{'}^{'''}$" as its only free variable; and   
    ($\beta$) roughly speaking: under the intended interpretation, $\phi_n$
    \emph{is} satisfied at all (denumerably) infinite sequences
    $f$ of classes of classes of individuals (in the domain of interpretation)
    whenever the first element of the sequence is the class of classes of such individuals that have
    exactly $n$ members.
  \item $\phi_n$ is an $\mathcal{L}_{GTC}$ sentence that is true
    under the intended interpretation, or, equivalently, is in $Tr$;
  \end{enumerate}
\end{enumerate}
To establish that an arithmetic definition of the required
(meta)metatheoretical sentential function $\psi(n)$ is indeed
available in ${T_M}_{GTC}$, let us firstly
assume that ${T_M}_{GTC}$ includes proper axioms
defining an order four constant, $Tr^{4}$, such that we may formally prove in ${T_M}_{GTC}$
the type-three variable equal to a given natural number $k$ in some context is
(or is not) the G\"odel number of an $\mathcal{L}_{GTC}$ sentence that is in $Tr$. To
introduce these hypothesised proper axioms, let's adopt some additional notation
for describing the expressions involved.
\begin{Definition}[$\epsilon_{Tr,i}^3$]
  For $i \in \mathbb{N}$, let $\epsilon_{Tr,i}^3$ be a (metametalinguistic) abbreviation for the
  metalinguistic expression $\ulcorner Tr^{4}(X^3_{i}) \urcorner$, noting that the
  "$i$" inside the corners is of course a metametalinguistic abbreviation for $i$ strokes.
\end{Definition}
Let us use $\phi_{Tr}(k)$ as the name of the
the following metametalinguistic sentential function:
\begin{Definition}[$\phi_{Tr}(k)$]
  Where $\phi_k$ is the $k^{\text{th}}$ expression in the above enumeration of
  $\mathcal{L}_{GTC}$ expressions, let $\phi_{Tr}(k)$ be
  the metametalinguistic sentential function that, for a definite choice of natural number for $k$,
  takes as value an ${\mathcal{L}_M}_{GTC}$ expression defined as follows:
  \begin{equation}
    \textstyle \phi_{Tr}(k) =_{df} \ulcorner \bigcup_{1}^{3}(\iota_{k}\ .\
    (\phi_k\ \in\ Tr\ .\ \epsilon_{Tr,1}^3) + (\phi_k\ \notin\ Tr\ .\ \overline{\epsilon_{Tr,1}^3})) \urcorner 
  \end{equation}
\end{Definition}
For arbitrary natural number $k$, the $k^{\text{th}}$ instance of the proper axiom
defining $Tr^4$ is as follows:
\begin{flalign*}
  (Tr^4_{k})\quad \phi_{Tr}(k)
\end{flalign*}

To further consider an approach to defining  $\psi(n)$,
let us also define the following metametalinguistic sentential function
($\psi^{'}(k)$) that, for arbitrary natural number $k$, takes as value
an ${\mathcal{L}_M}_{GTC}$ sentence:
\begin{align*}
  \textstyle \psi^{'}(k) =_{df} \ulcorner \bigcup_{1}^{3}(\iota_{k}\ .\ \epsilon_{Tr,1}^3) \urcorner
\end{align*}
With $\psi^{'}(k)$) thus defined, the remaining details
for the definition of the required function $\psi(n)$ concern
only routine matters of the arithmetization of the syntax of $\mathcal{L}_{GTC}$ -
such as the arithmetic definition of the class of G\"odel numbers of the class of
$\mathcal{L}_{GTC}$ sentential functions, the definition of the arithmetic
relation that corresponds to the metamathematical relation that holds
when $x$ is an $\mathcal{L}_{GTC}$ free variable in the $\mathcal{L}_{GTC}$
sentential function $y$, the arithmetization of an $\mathcal{L}_{GTC}$
expression $x$ that is the negation of another  $\mathcal{L}_{GTC}$
expression $y$ and so on. Methods for addressing these problems are all well
covered in \cite{godel1931}.

Let us grant then Tarski's assertion that, given
the assumption that a well defined class $Tr$ exists, and
the other syntactical assumptions set out above, we may
infer that there exists the required (arithmetically defined)
sentential function "$\psi(n)$" - that is truth functionally equivalent
to ``$\bigcup_{1}^{3}(\iota_{n} . \phi_n)\ \overline{\in}\ Tr$''. Rejoining
Tarski's demonstration at this point, we have:
\begin{quote}
  (1) \indent \emph{for any} $n$,  $\bigcup_{1}^{3}(\iota_{n} . \phi_n)\ \overline{\in}\ Tr$ \emph{if and only if}  $\psi(n)$ (\cite{tarski1936}: 250).
\end{quote}
Since $\psi(n)$ will be one of the expressions in the aforementioned enumeration
of the expression of the object language ($\phi_0$, $\phi_1$, \ldots, $\phi_n$, \ldots),
such as the $k^{\text{th}}$ ($\psi(n) = \phi_k$); if we substitute $k$ for $n$ in $(1)$ we obtain:
\begin{quote}
  (2) \indent $\bigcup_{1}^{3}(\iota_{k} . \phi_k)\ \overline{\in}\ Tr$ \emph{if and only if}  $\psi(k)$ (\cite{tarski1936}: 250).
\end{quote}
Since $\bigcup_{1}^{3}(\iota_{k} . \phi_k)$, for a definite choice of $k$,
is a metalinguistic name of a sentence of the object language, we may
also obtain the following as an instance of the $\alpha$ sentences
identified in Convention T (by choosing $\bigcup_{1}^{3}(\iota_{k} . \phi_k)$ for $x$ and $\psi(k)$ for $p$):
\begin{quote}
  (3) \indent $\bigcup_{1}^{3}(\iota_{k} . \phi_k)\ \in\ Tr$ \emph{if and only if}  $\psi(k)$ (\cite{tarski1936}: 250).
\end{quote}
To summarise the result obtained, and provide an alternative statement of the Theorem established, the following definition may be used:
\begin{Definition}\label{definition_liar_sentence}
  [Liar Sentence] An object language ($\mathcal{L}$) sentence $\mathcal{B}$ shall be referred to as a liar sentence
  for a formal theory ($T_{\mathcal{L}}$) in this language under an interpretation $\mathfrak{M}$, iff
  (with semantic notions defined as per Tarski's theory):
  the hypothesis that $\mathcal{B}$ is true under $\mathfrak{M}$ ($\vDash_{\mathfrak{M}} \mathcal{B}$) implies,
  in a formal metatheory $M$ in which the interpretation is defined, 
  that $\mathcal{B}$ is not true under $\mathfrak{M}$ ($\nvDash_{\mathfrak{M}} \mathcal{B}$) and vice versa:
  \begin{equation}
    (\vDash_{\mathfrak{M}} \mathcal{B}) \leftrightarrow (\nvDash_{\mathfrak{M}} \mathcal{B})
  \end{equation}
\end{Definition}
Thus, the sentence named by Equation (3) is a liar sentence for our formal metatheory $M$ of $T_{GTC}$
under the hypothesised interpretation, since, in our metametatheory: (2), (3),
and an instance of the tautology $(A \sim B) \land (C \sim B) \to (A \sim C)$ yield (via
modus ponens):
\begin{equation}
  \textstyle \bigcup_{1}^{3}(\iota_{k} . \phi_k)\ \in\ Tr\ \text{iff}\ \bigcup_{1}^{3}(\iota_{k} . \phi_k)\ \overline{\in}\ Tr
\end{equation}
For future reference we may note an alternative formulation of Tarski's (\cite{tarski1936}) Theorem I (\S 5):
\begin{Proposition}\label{proposition_liar} For any formal theory $T$ in scope for Theorem I (\S 5),
  the hypothesis that there exists a well defined class $Tr$ of $T$ sentences that are true under the intended interpretation
  implies that there exists a sentence $\mathcal{B}$ in the language of $T$ that is a liar sentence,
  in the sense of Definition \ref{definition_liar_sentence}, for the theory $T$.
\end{Proposition}

Subsequent to proving the above Theorem, Tarski presented a revised Theorem discussed in the following subsection.

\subsection{Tarski's revised Theorem on languages of infinite order}\label{section_languages_infinite_revised}

Reflecting on the possibility of introducing transfinite variables into the metalanguage,
Tarski concluded that:
\begin{quote}
  \ldots it is always possible to con-struct the metalanguage in such a way that it contains variables
  of higher order than all the variables of the language studied. (\cite{tarski1936}: 271-2).
\end{quote}
In light of which, Tarski retracted the above Theorem and replaced it with the following:
\begin{quote}
  A. For every formalized language a formally correct and materi-ally adequate definition of
  true sentence can be constructed in the metalanguage with the help only of general logical expressions,
  of expressions of the language itself, and of terms from the morphology of language -
  but under the condition that the metalanguage possesses
  a higher order than the language which is the object of investigation.\newline
  B. If the order of the metalanguage is at most equal to that of the language itself, such a definition cannot be constructed.
  (\cite{tarski1936}: 273).
\end{quote}

In considering Tarski's revised position let's focus on the case of the general theory of classes discussed above.
The proposed metalanguage is to be comprised of - in addition to the object language $O$, the general theory of classes,
and the appropriate special symbols and proper axioms for describing $O$ (\cite{tarski1936}: 170-206 etc.) -
transfinite variables. Since there are at most a denumerable number of true $T_{GTC}$ formulae, let us assume
that the hypothesised metalanguage contains, in addition to a denumerable supply of variables
of each finite order, a denumerable supply of variables of order $\omega$, say $X^{\omega}_{'}$, $X^{\omega}_{''}$ etc.
As Tarski notes  some sentence forming functors of order
$\omega+1$ will be required (\cite{tarski1936}: 270) to formulate statements
involving the new variables. Tarski argues that with the introduction of the variables of transfinite
order it also becomes necessary to introduce variables and (sentence) functors of
indefinite order:
\begin{quote}
  We must introduce into the
  languages variables of indefinite order which, so to speak, 'run
  through' all possible orders, which can occur as functors or
  arguments in sentential functions without regard to the order of
  the remaining signs occurring in these functions, and which at
  the same time may be both functors and arguments in the same
  sentential functions. With such variables we must proceed with
  the greatest caution if we are not to become entangled in
  antinomies like the famous antinomy of the class of all classes
  which are not members of themselves. Special care must be
  taken in formulating the rule of substitution for languages which
  contain such variables and in describing the axioms which we
  have called pseudodefinitions.
  (\cite{tarski1936}: 271).
\end{quote}

To provide some context for this remark,
it may be useful to consider the stages view of the formation of classes (\cite{shoenfield1967}: \S9.1).
Although this view is more often used to develop an account of ZFC, we may also use it
as follows to develop an account of the General Theory of Classes under consideration.

Consider firstly the following fact about the
General Theory of Classes:
\begin{Fact}\label{fact_gtc_grammar}
  Under the formation rules of the General Theory of Classes,
  for any variables $X$ and $Y$, of finite order, the atomic
  formula $\ulcorner X Y \urcorner$
  is well formed if and only if (for some natural number
  $n$): $Y$ is of order $n$ and $X$ is of order $n+1$ (\cite{tarski1936}: 242).
\end{Fact}
This fact may be explained in terms of the stages
view of class formation as follows:
\begin{Hypothesis}\label{hypothesis_stages}
  If we assume that under the intended
  interpretation of the General Theory of Classes
  (for some natural number
  $n$): the $n^{th}$ order variables range over
  the collections of objects formed as sets
  at the $n^{th}$ stage, then the formation rules
  of this theory (Fact \ref{fact_gtc_grammar}),
  imply that a collection of objects formed as a set
  at the $(n+1)^{th}$ stage is only gathered together
  as a single object that forms a set
  \emph{after} the gathering of objects into sets
  at the $(n)^{th}$ stage is complete for \emph{all} of the
  latter collections.
\end{Hypothesis}
The following example may serve to illustrate
this view of class formation:
\begin{Example}\label{example_stages}
  According to the stages view (Hypothesis \ref{hypothesis_stages}),
  the classes of classes of individuals - which,
  under the intended interpretation, are ranged over
  by the second order variables ($X_{'}^{''}$, $X_{''}^{''}$, \ldots) -
  are only gathered together
  as single objects that form sets \emph{after}
  all of the classes of individuals - which,
  under the intended interpretation, are ranged over
  by the first order variables ($X_{'}^{''}$, $X_{''}^{''}$, \ldots) -
  are formed.
\end{Example}
By this account, heterogeneous classes -
that include sets as members that are ranged over by
variables of different orders - \emph{never} come into existence
because such a possibility is excluded by the assumption
(Hypothesis \ref{hypothesis_stages}) that classes are formed in stages.

To illustrate with an example again:
\begin{Example}\label{example_stages_conflict}
  According to the stages view (Hypothesis \ref{hypothesis_stages}),
  there is no class that includes as members both an individual $a$
  and the singleton class consisting solely of
  this individual ($\{a \}$):
  \begin{equation*}
    (\text{Prohibited by Hypothesis \ref{hypothesis_stages}})\quad  a \in X \land \{a \} \in X
  \end{equation*}
  The existence of such a
  set would require the gathering of
  a collection of a \emph{class of classes of individuals} as a single
  entity (a set), \emph{prior} to
  the completion of the formation of \emph{all} classes of
  individuals (since the hypothesised class
  $X$ would itself be one of these very same sets being formed,
  albeit, per impossible, a set that mingles individuals with
  classes of individuals).
\end{Example}

In certain cases, classes and relations
which are heterogeneous in informal usage, can
be simulated by a suitably defined homogeneous class or relation
as G\"odel points out:
\begin{quote}
  Nonhomogeneous relations, too, can be defined in this manner; for example, a
  relation between individuals and classes can be defined to be a class of elements of
  the form $((x_2), ((x_1), x_2))$ (\cite{godel1931}: 153, fn 18).
\end{quote}
G\"odel, using here the symbolism for his theory $P$,
is explaining the possibility of defining within $P$
a relation that corresponds to an \emph{informal},
heterogeneous relation that holds between
individuals and classes of individuals, if for example
we assert that a certain number, such as two - viewed as an individual -
is a member of a class of these individuals (such as the class of prime
numbers). The relation ($((x_2), ((x_1), x_2))$) thus defined within
the system $P$, which corresponds, in terms of content, to this informal relation,
is of course a \emph{homogeneous} class as defined in terms of the
symbolism of system $P$.

In other cases however it may not be
possible to construct a homogeneous formal
definition for a class or relation which is
heterogeneous in informal usage.  von Neumann's definition of the
ordinals (\cite{shoenfield1967}: \S9.3) can't be formulated
in the general theory of classes presented above without some modifications
along the lines of Tarski's suggestion. On the von Neumann approach,
the ordinals are built, commencing with the null class, by forming the successor
($S$) of any given ordinal $\alpha$ as the union of the ordinal
and the singleton set of this ordinal itself:
\begin{equation*}
  (\text{von Neumann successor for ordinals})\quad  S(\alpha) = \alpha \cup \{\alpha\}
\end{equation*}
As Example \ref{example_stages} explains, this
definition of the successor operation for ordinals is
not compatible with the above General Theory of Classes. Similarly,
on the von Neumann approach, the least transfinite ordinal $\omega$ is defined as follows:
\begin{equation*}
  (\text{von Neumann definition of $\omega$})\quad \emptyset \in \omega \land (\forall \sigma)(\sigma \in \omega \rightarrow S(\sigma) \in \omega)
\end{equation*}
For the same reason as set out in Example \ref{example_stages_conflict},
there can be no such set as $\omega$ defined in the above
General Theory of Classes. The above definition of $\omega$ also illustrates a related
complication: as $\omega$ is, by hypothesis, a limit ordinal - it has no immediate predecessor -
the pattern of the existing hierarchy of stages cannot be relied upon
to determine the order of $\omega$'s members; the ordinals prior to
$\omega$, being finite, all have immediate predecessors.

In the light of these considerations, let us accept Tarski's recommendation and
adopt the following hypothesis with respect to the metatheory:
\begin{Hypothesis}\label{hypothesis_ordinal_variables}
  A denumerable supply of variables ($Y_{'}$, $Y_{''}$ \ldots) that
  (under the intended interpretation) range over the ordinals
  (suitably defined), together
  with appropriate axioms for the use of these symbols, have been added to
  our metatheory $M$ of the general theory of classes without introducing any
  pathology into the system.
\end{Hypothesis}
Clearly the hypothesised axioms
need to exclude various unwarranted substitutions and inferences - since under the intended interpretation,
the new variables that range over the finite ordinals will vary over all orders of the classes ranged
over by the existing variables of the object
language ($X^{'}_{i}$, $X^{''}_{j}$, $X^{'''}_{k}$ etc.).\footnote{To keep matters simple,
I will not, as Tarski proposed, re-index the variables so that variables for individuals are order $0$ ($X^{0}_{i}$) etc.}

For a specific system that illustrates this idea we might consider the set theory Quine
presents in \emph{Mathematical Logic} (\cite{quine1981}): the stratification requirement
imposed on the axioms of members (\cite{quine1981}: \S 28) closely parallels the syntax of the
atomic formulae of the General Theory of Classes; the absence of a syntactically imposed
typing of the provides the variables in other contexts with the "indefinite order" character, so that
von Neumann's ordinals may be defined in the system (\cite{quine1981}: \S 45). Of course, Quine's
is a first-order theory, so the idea that it illustrates an approach to defining the sought
after metatheory should not be taken too literally.

Let us assume then that we have added to $M$ proper axioms allowing for use of
these new (order-free) ordinal variables and an abbreviation for the successor function for ordinals.

In summary, one approach to making precise Tarski's revised position on the definability of truth
for the object language of interest, the general theory of classes, is this:
\begin{enumerate}
\item The metalanguage $M$ includes newly added order-free variables, and
  the newly added transfinite $M$ symbol for $Tr$ ($Tr^\omega$ to avoid confusion),
  together with appropriate axioms defining proper use of these symbols;
\item Under the intended interpretation of $M$, $Tr^\omega$ names the transfinite ordinal
  that contains all, and only, the members of a denumerable sequence
  of ordinals ($\sigma_0$, \ldots $\sigma_n$, \ldots) that
  name the object language sentences that are true
  under the intended interpretation; and
\item A formally correct and materially adequate definition of the true sentences
  of the general theory of classes can be constructed in the (transfinite) metalanguage $M$ thus loosely described. 
\end{enumerate}

Let us now consider the question of whether this approach to defining truth for the general theory of classes avoids
the liar paradox constructed in \S \ref{section_languages_infinite_paradox}.
\begin{Proposition}\label{proposition_tarski_revised} Tarski's revised approach to defining
  truth for the general theory of classes (\cite{tarski1936}: \S7) fails to avoid the paradox
  identified in his initial consideration  (\cite{tarski1936}: \S5) of this problem.
\end{Proposition}
\begin{proof}
  For brevity, I assume as proven
  with respect to our metatheory $M$ the following result
  (for which Shoenfield presents a proof with respect to ZFC):
  \begin{Theorem}\label{theorem_zermelo}
    If we assume the axiom of choice in $M$,
    then we may prove Zermelo's well-ordering theorem:
    \begin{quote}
      For every set $x$, there is a bijective mapping from an ordinal to
      $x$ (\cite{shoenfield1967}: 253).
    \end{quote}
  \end{Theorem}
  Since there are at most a denumerable number of ordinals in $Tr^\omega$, the
  set is not only, by the above theorem, well-ordered but we may use
  this ordering to map the members of $Tr^\omega$ to the finite ordinals.
  Hence, hereon, to reduce circumlocutions, I will, wherever convenient,
  treat $Tr^\omega$ itself as a set of finite ordinals (rather
  than speaking of the set of finite ordinals that corresponds, in the
  manner just indicated, to $Tr^\omega$).
  
  If we thus take all of the members of $Tr^\omega$ to be finite ordinals, Theorem \ref{theorem_zermelo}
  implies that there exists a mapping of all the these to the corresponding
  finite classes (of all classes of $k$ individuals, or $l$ individuals etc.)
  For brevity I assume that a definition of the
  relation of similarity (\cite{shoenfield1967}: 252), which may be used to
  establish this mapping, is available in $M$:
  \begin{Definition}[$\text{Sm}(y, x)$]
    Where $y$ and $x$ are $M$ variables for a finite ordinal and
    an order three class (a class of classes of individuals),
    e.g.  $Y_{1}$ and $X^3_{1}$,
    let $\text{Sm}(y, x)$ be the $M$ sentential function (relation sign) which,
    under the intended interpretation, asserts the existence of
    a bijective mapping between $y$ and $x$  
  \end{Definition}
  I assume throughout the following that some fixed
  G\"odel numbering of the object system, the general theory of classes,
  is given:
  \begin{Definition}[$\psi_k$]
    For the given G\"odel numbering of the object system $O$
    and arbitrary natural number $k$,
    $\psi_k$ shall be the $k^{th}$ wfs of $O$ according to this
    enumeration.\footnote{Quine's corners would have been appropriate
    to make this definition more precise, where it not for the risk that some readers
    would mistake the corners for a function that returns the G\"odel number
    of the enclosed expression.}
  \end{Definition}

  Hereon where convenient, to improve the readability of $M$ formulae, I will use various
  abbreviations for $M$ expressions in the customary manner - introducing symbols for
  the conjunction, the conditional, the biconditional and existential quantification
  (".", "$\rightarrow$", "$\leftrightarrow$", "$\exists$").
  
  A proper axiom schema for $M$ that, under the intended interpretation,
  for arbitrary choice of natural number $k$,
  defines whether the relevant finite ordinal is or is not a member of $Tr^\omega$ (where $\sigma$ is one of the
  ordinal variables introduced into $M$ as per Hypothesis \ref{hypothesis_ordinal_variables}):
  \begin{enumerate}
  \item  If $k$ is the G\"odel number of an $O$ sentence, the axiom schema defining $Tr^{\omega}$
    for such choices of $k$ is as follows:
    \begin{flalign*}
      (Tr^{\omega}_k)\quad \prod\sigma \{[\sigma \in Tr^\omega] \leftrightarrow 
      (\exists X^3_{j}) [X^3_{j} = k\ .\ \text{Sm}(\sigma, X^3_{j})\ .\ \psi_k]\} &&
    \end{flalign*}
  \item If $k$ is not the G\"odel number of an $O$ sentence, the axiom schema defining $Tr^{\omega}$
    for such choices of $k$ is as follows:
    \begin{flalign*}
      (Tr^{\omega}_k)\quad \prod\sigma \{ (\exists X^3_{j}) [X^3_{j} = k\ .\ \text{Sm}(\sigma, X^3_{j})\ \rightarrow N [\sigma \in Tr^\omega]]\} &&
    \end{flalign*}
  \end{enumerate}
  If the class defined by $Tr^\omega$ under the intended interpretation
  of the metatheory $M$ is well-defined, then, there exists
  an instance of an $M$ comprehension axiom asserting
  the existence of a class (of classes of classes of individuals) that contains all,
  and only, the finite classes of $k_0$, $k_1$, \ldots individuals -
  with each $k_i$ being the G\"odel number of a true sentence - that are similar
  to the finite ordinals that are members of $Tr^{\omega}$. Where convenient hereon I will
  use G\"odel's statement of the comprehension axioms (\cite{godel1931}: 155, Ax IV.1).
  Whilst Tarski Psudodefinition axioms (\S \ref{section_languages_infinite_paradox})
  are similar, G\"odel's format is more convenient for our purpose. In somewhat modernised notation,
  the required instance of these axioms is as follows.\footnote{$X^4_{i}$ must not
  occur free to the right of the biconditional; to obtain an instance of the axiom, abbreviations must
  be eliminated and all official parentheses restored (\cite{godel1931}: 155).}
  \begin{flalign*}
    (Tr^4)\quad (\exists X^4_{i})\{(\forall X^3_{j}) [(X^3_{j} \in X^4_{i}) & \leftrightarrow \\
      & (\exists \sigma) (\sigma \in Tr^{\omega}\ .\ \text{Sm}(\sigma, X^3_{j}))]\}  &&
  \end{flalign*}

  As $Tr^4$  is an $M$ Theorem,
  this implies the following metatheorem concerning
  a system $M'$ related to $M$ (cf. \cite{tarski1936}:256 Theorem III):
  \begin{quote}
    \textbf{Metatheorem $III^{\prime}$} 
    Let $M'$ be the metatheory that results when $M$ is modified by addition of (i) a new primitive symbol for an order four sign, $Tr^4$,
    and (ii) proper axioms corresponding to $T$ biconditionals for the object theory $O$ (the general theory of classes),
    obtained from the following axiom schemata:
    \begin{enumerate}
    \item  If $k$ is the G\"odel number of an $O$ wff, the axiom schema defining $Tr^{4}$
      for such choices of $k$ is as follows ():
      \begin{equation*}
        (Tr^4 Ax_k)\quad (\exists X^3_{i})(X^3_{i} = k\ .\ X^3_{i} \in Tr^4) \leftrightarrow \psi_k
      \end{equation*}
    \item If $k$ is not the G\"odel number of an $O$ wff, the axiom schema defining $Tr^{4}$
      for such choices of $k$ is as follows:
      \begin{equation*}
        (Tr^4 Ax_k)\quad \prod X^3_{i} [(X^3_{i} = k)\ \rightarrow N (X^3_{i} \in Tr^{4})]
      \end{equation*}
    \end{enumerate}

    If the class defined by $Tr^\omega$ under the intended interpretation
    of the metatheory $M$ is well-defined, then so too is the class defined by $Tr^4$ under the intended interpretation
    of the metatheory $M'$ 
  \end{quote}
  If however we now compare the proper axioms added as instances of $Tr^4 Ax_k$ with the proof of
  \cite{tarski1936}: \S5 Theorem I, we see that Theorem $III^{\prime}$ shows that the paradox that established
  Tarski's original Undefinability Theorem may easily be reproduced
  with the aid of the newly added axioms. For a more modern style of proof, we could add numerals for natural numbers,
  together with proper axioms defining these, a sentential function $\mathcal{D}(x, y)$ representing the diagonal function.
  Let $m$ be the G\"odel number of the following sentential function:
  \begin{equation*}
    (\mathcal{M}(X^3_{1}))\quad  (\forall X^3_{2}) \{\mathcal{D}(X^3_{1}, X^3_{2}) \rightarrow [X^3_{2} \notin Tr^4]\}
  \end{equation*}
  Let $q$ be the G\"odel number of the sentence results when the numeral for $m$ is substituted for all
  free occurrences of $X^3_{1}$ in $(\mathcal{M}(X^3_{1}))$. Then we may prove in $M$ that,
  under the intended interpretation, the following sentence is a liar sentence
  for this system:
  \begin{equation*}
    (\forall X^3_{2}) \{\mathcal{D}(\overline{m}, X^3_{2}) \rightarrow [X^3_{2} \notin Tr^4]\} \leftrightarrow (\overline{q} \notin Tr^4)
  \end{equation*}
\end{proof}
As an easy corollary we have:
\begin{Corollary}
  If the above proof is examined it may be confirmed that, as a rule, it applies to
  any classical set theory in scope for both G\"odel's proof and Tarski's theory of truth
  for formal languages.
\end{Corollary}

At this point we can consider application of Tarski's revised Undefinability Theorem to the case of first-order arithmetic.
\section{Tarski's Theorem and first-order arithmetic}
For the case of first-order arithmetic the
slogan version of Tarski's revised Undefinability Theorem is this:
``arithmetic truth is not arithmetically definable'' (\cite{mendelson2015}: 220, cf. \cite{boolosburgessjeffrey2007}: 222). That is,
there cannot be any formula $T(x)$ of
first-order arithmetic such that, under
the intended interpretation, $\ulcorner T(\overline{n}) \urcorner$ is true if and only if $\overline{n}$ is the numeral for the G\"odel number of a sentence of first-order arithmetic
that is true under this interpretation.
\par
Since Mendelson's \cite{mendelson2015} treatment of this material is first rate, including very precise definitions of the required content,
I shall take this as a definitive statement of
the application of Tarski's Theorem to first-order arithmetic. Mendelson formulates Tarski's Theorem for first-order arithmetic as follows:
\begin{quote}
  \textbf{Corollary 3.44 [Tarski's Theorem (1936)]} Let $Tr$ be the set of G\"odel numbers of wfs [well-formed formulas] of $S$ [first-order arithmetic]
  that are true for the standard interpretation. Then $Tr$ is not arithmetical. (\cite{mendelson2015}: 220 modified through interpolation in square brackets.)
\end{quote}
For reference, Mendelson's surprisingly succinct proof of Tarski's Theorem for first-order arithmetic (Corollary 3.44) is this:\footnote{The reader who is unfamiliar with \cite{mendelson2015} will find a summary of the notation / terminology used in the proof, including a statement of
the Proposition 3.43 mentioned in the proof, in \S \ref{section_notation_terminology_corollary}}
\begin{quote}
  Let $\mathcal{N}$ be the extension of $S$ that has as proper axioms all those wfs that are true for the standard interpretation.
  Since every theorem of $\mathcal{N}$ must be true for the standard interpretation, the theorems of $\mathcal{N}$ are identical with the axioms of $\mathcal{N}$.
  Hence $T_{\mathcal{N}} = Tr$. Thus, for any closed wf $\mathcal{B}$, $\mathcal{B}$ holds for the standard interpretation if and only if
  $\vdash_{\mathcal{N}}\mathcal{B}$. It follows that a set $B$ is arithmetical if and only if the property $x \in B$ is expressible in $\mathcal{N}$. We may assume that
  $\mathcal{N}$ is consistent because it has the standard interpretation as a model. Since every recursive function is representable in $S$, every recursive function is representable
  in $\mathcal{N}$ and, therefore, $D$ is representable in $\mathcal{N}$. By Proposition 3.43 $x \in Tr$ is not expressible in $\mathcal{N}$. Hence $Tr$ is not arithmetical.
  (This result can be roughly paraphrased by saying that the notion of arithmetical truth is not arithmetically definable.)
  (\cite{mendelson2015}: 220).
\end{quote}
If we compare Corollary 3.44 with Tarski's \S5 Theorem I reproduced
above it is apparent that a rather interesting transformation of the Theorem has occurred.
In the original statement of the Theorem, it is the hypothesis that $Tr$ is defined that is to be reduced to absurdity, whereas for the case of first-order
arithmetic (Corollary 3.44) it is instead the hypothesis that $Tr$ is arithmetical that is to be reduced to absurdity. In view of the contents of
Proposition \ref{proposition_liar} this raises an interesting issue about the proof of Corollary 3.44. To examine this issue, let's consider
the hypothesis that Proposition \ref{proposition_liar}, with appropriate changes, applies to first-order arithmetic:
\begin{Hypothesis}\label{hypothesis_first_order_arithmetic} The hypothesis that $Tr$ is well defined
  for the language of first-order arithmetic implies that the class $Tr$ is arithmetical and
  hence that a liar sentence (Definition \ref{definition_liar_sentence}) may be defined in this theory.
\end{Hypothesis}
Hypothesis \ref{hypothesis_first_order_arithmetic} is introduced for consideration to highlight
a logical flaw in the proof of Corollary 3.44 that may be stated as follows:
\begin{Proposition}\label{proposition_mendelson_tarski}
  The above proof of Tarski's Theorem for first-order arithmetic (Corollary 3.44)
  \emph{assumes as a premise} that Hypothesis \ref{hypothesis_first_order_arithmetic} is false.
\end{Proposition}
\begin{proof}
  The proof by contradiction is as follows. Assume in order to derive a contradiction that the \emph{falsity} of Hypothesis \ref{hypothesis_first_order_arithmetic} 
  is \emph{not} necessary
  in order for the above proof of Tarski's Theorem for first-order arithmetic (Corollary 3.44) to be valid. Thus, we may assume that Hypothesis \ref{hypothesis_first_order_arithmetic} is
  true and yet the above proof of Corollary 3.44 yields the desired conclusion from the indicated premises. If however we examine the premises of the proof we may confirm that
  among these is the assumption that the set $Tr$ is well defined. (The assumption is evident in the assumption that the definition of $\mathcal{N}$ is well defined -
  so that the set of $\mathcal{N}$ theorems etc is well defined.)  This assumption, in conjunction with the assumption that Hypothesis \ref{hypothesis_first_order_arithmetic} is true,
  yields via modus ponens the conclusion that the set $Tr$ is arithmetically defined -
  a contradiction to the stated conclusion of the proof: that  $x \in Tr$ is not expressible in $\mathcal{N}$.
\end{proof}
Proposition \ref{proposition_mendelson_tarski} yields the following corollary:
\begin{Corollary}\label{corollary_mendelson_petitio}
  The above proof of Tarski's Theorem for first-order arithmetic (Corollary 3.44)
  fails on the grounds that the result that is to be established is assumed as a premise.
\end{Corollary}
\begin{proof}
  Let $A$ and $\lnot B$ be the following propositions:
  \begin{description}
  \item[$A$] The set $Tr$ is well defined.
  \item[$\lnot B$] The set $Tr$ is \emph{not} arithmetically defined.
  \end{description}
  Clearly, the proof of Tarski's Theorem for first-order arithmetic (Corollary 3.44) aims to establish the truth of $\lnot B$.
  That $A$ is essential for the proof is confirmed by inspection of, for example, the definition of the formal system $\mathcal{N}$ - if
  the set $Tr$ is \emph{not} well defined then neither is the formal system  $\mathcal{N}$ and the proof collapses. Hence, if the proof is sound $A$ is true
  and thus, by conjunction introduction, the intended conclusion of the proof $\lnot B$ is logically equivalent to the conjunction $A \land \lnot B$.

  The proposition that this conjunction ($A \land \lnot B$) is true is also equivalent however to the proposition that Hypothesis \ref{hypothesis_first_order_arithmetic} ($\alpha$), $A \Rightarrow B$,
  is false:
  \begin{flalign*}
    (A \land \lnot B) \Leftrightarrow \lnot (A \Rightarrow B)  &&
  \end{flalign*}
  By Proposition \ref{proposition_mendelson_tarski}, the proof of Corollary 3.44 \emph{assumes}
  that Hypothesis \ref{hypothesis_first_order_arithmetic} is false. Hence, the proof of Corollary 3.44 assumes
  as a premise a hypothesis that is equivalent to the result that is to be proven.
\end{proof}
\section{Conclusion}
The above presentation would probably be briefer and easier to follow if,
rather than following Tarski's original proof of the Undefinability Theorem (\cite{tarski1936} \S 5 Theorem I),
a more contemporary style was adopted. I chose the former course to emphasise that much of the
detail of the argument presented is old; my initial rendering of these
claims should have indicated that Tarski had rescinded this position.

For most contemporary readers, results concerning the status of first-order theories, and of
$ZFC$ and $PA$ in particular, are
of much greater interest than those concerning higher-order theories. Higher-order theories do however
continue to play an important role in some areas, so the above-presented claims about these
systems will hopefully be found useful by readers with such interests.
\appendix
\section{Notation and terminology used in the proof of Corollary 3.44}\label{section_notation_terminology_corollary}
The reader who is unfamiliar with Mendelson's excellent publication
may wish to note the following points concerning the notation / terminology used in the proof of Corollary 3.44:
\begin{enumerate}
\item $S$ is Mendelson's version of first-order number theory. The proof only uses properties of $S$ that are shared by any standard version of first-order number theory,
  but for the curious: Mendelson's first-order theories (\cite{mendelson2015}: \S 2.3) are defined with the aid of five axiom schema and two rules of inference, modus ponens
  and generalisation. As proper axioms for first order arithmetic  (\cite{mendelson2015}: \S 3.1), eight open formulae and an induction schema are used.
\item The language of arithmetic $\mathcal{L}_A$ (\cite{mendelson2015}: \S 3.1) is taken to include a single predicate letter for equality (``=''), a single constant for zero (``0''), and three function
  letters for the successor function (``$\ '\ $''), addition (``+'') and multiplication (``.''), together with a denumerable supply of individual variables ($x_{1}, \ldots, x_{n}$)
  and symbols for negation (``$\lnot$'') and the conditional (``$\Rightarrow$'').\footnote{For some symbols, only the unofficial versions are listed since the official versions are not used in the proof.}
\item The formal theory $\mathcal{N}$: The proof of Corollary 3.44 does not focus directly on $S$ but a formal theory $\mathcal{N}$, defined in the proof as the following extension of $S$:
  \begin{quote}
    Let $\mathcal{N}$ be the extension of $S$ that has as proper axioms all those wfs [well-formed formulae] that are true for the standard interpretation.
    (\cite{mendelson2015}: 220)
  \end{quote}
\item A standard definition of primitive recursive and recursive number-theoretic functions and relations is used (\cite{mendelson2015}: \S 3.3).
\item The arithmetization of the formal theories $S$ / $\mathcal{N}$: the function $g$, defined as follows, is used to map the symbols, expressions and
  sequences of expressions of $S$ / $\mathcal{N}$ into
  the natural numbers  (\cite{mendelson2015}: 192-3):
  \begin{enumerate}
  \item (Symbols) $g(() = 3$, $g()) = 5$, $g(,) = 7$, $g(\lnot) = 9$, $g(\Rightarrow) = 11$, $g(\forall) = 11$, ``0'': $g(a_{1}) = 15$,
    ``$\prime$'': $g(f_{1}^{1}) = 49$, ``+'': $g(f_{1}^{2}) = 97$, ``$\cdot$'': $g(f_{2}^{2}) = 289$, ``='': $g(A_{1}^{2}) = 99$
  \item (Sequences of Symbols) The sequence of  $S$ / $\mathcal{N}$ symbols $u_0$ $u_1$ \ldots $u_j$ is mapped to the number $g(u_0 u_1\ldots u_j)$ defined as follows (where
    $p_{j}$ denotes the $j^{\text{th}}$ prime, with two being the zeroth prime):\\
    $g(u_0 u_1\ldots u_j) = 2^{g(u_{0})}\times 3^{g(u_{1})} \ldots \times p_{j}^{g(u_{j})}$
  \item (Sequences of Sequences of Symbols) If $e_0$, $e_1$, \ldots $e_j$ is a (finite) sequence of  expressions of $S$ / $\mathcal{N}$, this sequence is mapped to the number
    $g(e_{0}, e_{1},\ldots, e_{j})$ defined as follows (where
    $p_{j}$ denotes the $j^{\text{th}}$ prime, with two being the zeroth prime):\\
    $g(e_{0}, e_{1},\ldots, e_{j}) = 2^{g(e_{0})}\times 3^{g(e_{1})} \ldots \times p_{j}^{g(e_{j})}$
  \end{enumerate}
\item $T_{\mathcal{N}}$ is the set of G\"odel numbers of theorems of $\mathcal{N}$.
\item $\ulcorner \vdash_{\mathcal{N}} \mathcal{B} \urcorner$ denotes that $\mathcal{B}$ is a theorem of $\mathcal{N}$ (\cite{mendelson2015}: \S 1.4). (I use corners here and throughout for Quine's
  (\cite{quine1981}: \S6 ) quasi-quotation,
  whereas Mendelson uses ``$\ulcorner \mathcal{B} \urcorner$'' for the G\"odel number of the formula $\mathcal{B}$.)
\item Expressible: For a (first order) formal theory $K$ in the language of arithmetic $\mathcal{L}_A$
  \begin{quote}
    a number-theoretic relation $R$ of $n$ arguments  is
    \emph{expressible} in $K$ if and only if, there is a wf [well-formed formula] $\mathcal{B}(x_{1}, \ldots, x_{n})$ of $K$ with the free variables $x_{1}, \ldots, x_{n}$
    such that, for any natural numbers $k_{1}, \ldots, k_{n}$ the following hold:
    
    1. If $\mathcal{R}(k_{1}, \ldots, k_{n})$ is true then $\vdash_{K} \mathcal{B}(\overline{k}_{1}, \ldots, \overline{k}_{n})$.
    
    2. If $\mathcal{R}(k_{1}, \ldots, k_{n})$ is false then $\vdash_{K} \lnot \mathcal{B}(\overline{k}_{1}, \ldots, \overline{k}_{n})$.      
    (\cite{mendelson2015}: 169 modified through interpolation in square brackets.)
  \end{quote}
  where $\overline{k}_{n}$ is the $\mathcal{L}_A$ numeral for the number $k_{n}$, i.e. $0^{'\ldots k_{n}\text{times}}$ the $k_{n}^{\text{th}}$
  successor of zero.
\item Free for: Let $\mathcal{B}(x)$ be a wf in a first-order language with free occurrences of the variable $x$.
  With the notion of a term of a first-order language defined in the customary manner, a term $t$ is said to be ``free for $x$'' in $\mathcal{B}(x)$
  if no free occurrence of $x$ in $\mathcal{B}(x)$ lies within the scope of a quantifier $\forall{y}$ /  $\exists{y}$ for a variable $y$ that occurs in $t$.
  (\cite{mendelson2015}: 50)
\item Theory with equality: If $K$ is a first order theory with a symbol for equality (informally ``$=$''),
  \begin{quote}
    $K$ is called a \emph{first-order theory with equality} (or simply \emph{a theory with equality}) if the following are theorems of $K$:
    
    (A6) If $(\forall x_{1})x_{1} = x_{1}$ \indent (reflexivity of equality)
    
    (A7) $x = y \Rightarrow (\mathcal{B}(x,x) \Rightarrow \mathcal{B}(x,y)$. \indent (Substitutivity of equality)

    where $x$ and $y$ are variables, $\mathcal{B}(x,x)$ is any wf, and $\mathcal{B}(x,y)$ arises from
    $\mathcal{B}(x,x)$ by replacing some but not necessarily all, free occurrences of $x$ by $y$, with
    the proviso that $y$ is free for $x$ in $\mathcal{B}(x,x)$. Thus $\mathcal{B}(x,y)$ may or may not
    contain any free occurrences of $x$. (\cite{mendelson2015}: 93)
  \end{quote}
\item $(\exists_{1}{x})$: With $y$ the first new variable that does not occur in $B(x)$, $(\exists_{1}){x}B(x)$ is used for:
  \indent $(\exists{x})B(x) \land (\forall{x})(\forall{y})(B(x) \land B(y) \Rightarrow x = y)$ (\cite{mendelson2015}: 98)
\item Representable: For a first-order theory with equality $K$ in the language of arithmetic $\mathcal{L}_A$
  \begin{quote}
    A number-theoretic function $f$ of $n$ arguments  is
    said to be \emph{representable} in $K$ if and only if, there is a wf
    $\mathcal{B}(x_{1}, \ldots, x_{n}, y)$ of $K$ with the free variables $x_{1}, \ldots, x_{n}, y$
    such that, for any natural numbers $k_{1}, \ldots, k_{n}, m$ the following hold:
    
    1. If $f(k_{1}, \ldots, k_{n}) = m$ then $\vdash_{K} \mathcal{B}(\overline{k}_{1}, \ldots, \overline{k}_{n},  \overline{m})$.
    
    2. $\vdash_{K} (\exists_{1}{y}) \mathcal{B}(\overline{k}_{1}, \ldots, \overline{k}_{n}, y)$.      
    (\cite{mendelson2015}: 170 )
  \end{quote}
  where $\overline{k}_{n}$ is the $\mathcal{L}_A$ numeral for the number $k_{n}$, i.e. $0^{'\ldots k_{n}\text{times}}$ the $k_{n}^{\text{th}}$
  successor of zero.
\item The diagonal function $D$:
  \begin{quote}
    If $K$ is a theory in the language $\mathcal{L}_A$ \ldots the diagonal function $D$ has the property that,
    if $u$ is the G\"odel number of a wf $\mathcal{B}(x_{1})$, then $D(u)$ is the G\"odel number of the wf $\mathcal{B}(\overline{u})$
    (\cite{mendelson2015}: 205)
  \end{quote}
  ``$\mathcal{B}(x_{1})$'' here indicates that if $\mathcal{B}$ contains free occurrences of the variable $x_{1}$,
  $\ulcorner \mathcal{B}(\overline{u}) \urcorner $ is the result of substituting $\overline{u}$ for all free occurrences of $x_{1}$ in $\mathcal{B}$. $\mathcal{B}$ may contain
  other free variables. $\mathcal{B}$ may also contain no free occurrences of $x_{1}$ in which case ``$\mathcal{B}(\overline{u})$'' is simply $\mathcal{B}$ (\cite{mendelson2015}: 50).
\item \textbf{Proposition 3.43}
  ``Let $K$ be a consistent theory with equality in the language $\mathcal{L}_A$ in which the diagonal function $D$ is representable. then the property $x \in T_k$
  is not expressible in $K$.'' (\cite{mendelson2015}: 219)
\item  Arithmetical set: ``A set $B$ of natural numbers is said to be \emph{arithmetical} if there is a wf $\mathcal{B}(x)$ in the language $\mathcal{L}_A$,
  with one free variable $x$, such that, for every natural number $n$, $n \in B$ if and only if $\mathcal{B}(\overline{n})$ is true for the standard
  interpretation.'' (\cite{mendelson2015}: 219)
\end{enumerate}

\end{document}